\documentclass{amsart}
\usepackage{fullpage}
\usepackage{amsmath}
\usepackage{amssymb}
\begin{document}
\title[SECOND ORDER  NORMALIZATION ... -- P-R Drag]{SECOND ORDER  NORMALIZATION IN THE GENERALIZED PHOTOGRAVITATIONAL  RESTRICTED THREE BODY PROBLEM
  WITH POYNTING-ROBERTSON DRAG}
\author[B.S.KUSHVAH, J.P. SHARMA and
 B.ISHWAR]{ B.S.Kushvah$^1$, J.P. Sharma$^2$ and B.Ishwar$^3$ \\
1. JRF  DST Project, 2.  Co-P.I   DST Project  and  3. P.I. DST Project\\
    University Department of Mathematics,\\
    B.R.A. Bihar University Muzaffarpur-842001,
    Bihar (India)\\Email:ishwar\_ bhola@hotmail.com, bskush@hotmail.com}
\maketitle
\begin{abstract}
In this paper we have performed second order normalization in the generalized photogravitaional restricted three body problem with Poynting-Robertson drag. We have performed Birkhoff's normalization of the Hamiltonian. For this we have utilised Henrard's  method and expanded the coordinates of the third body in Double d'Alembert series. We have found the values of first and second order components. The second order components are obtained as solutions of the two partial differential equations. We have employed the first condition of KAM theorem in solving these equations. The first and second order components are affected by radiation pressure, oblateness and P-R drag. Finaly we obtained the third order part $H_3$ of the Hamiltonian in $I_1^{1/2}I_2^{1/2}$ zero.\end{abstract}
\section{Introduction}
The restricted three body problem describes the motion of an
infinitesimal mass moving under the gravitational effect of the two
finite masses, called primaries, which move in circular orbits
around their centre of mass on account of their mutual attraction
and the infinitesimal mass not influencing the motion of the
primaries. The classical restricted three body problem is
generalized to include the force of radiation pressure, the Poynting-Robertson effect and oblateness effect.

J. H. Poynting (1903) considered the effect of the absorption and
subsequent re-emission of sunlight by small isolated particles in
the solar system. His work was later modified by H. P. Robertson
(1937) who used a precise relativistic treatments of the first order
in the ratio of he velocity of the particle to that of light.

The effect of radiation pressure and P-R. drag in the restricted
three body problem has been studied by Colombo {\it et al.} (1966),
Chernikov Yu. A. (1970) and Schuerman (1980) who discussed the
position as well as the stability of the Lagrangian equilibrium
points when radiation pressure, P-R drag force are included. Murray
C. D. (1994) systematically discussed the dynamical effect of
general drag in the planar circular restricted three body problem,
Liou J. C. {\it et al.} (1995) examined the effect of radiation
pressure, P-R drag and solar wind drag in the restricted three body problem.

Moser's conditions (1962), Arnold's theorem (1961) and Liapunov's
theorem (1956) played a significant role in deciding the nonlinear
stability of an equilibrium point. Applying Arnold's theorem (1961),
Leontovic (1962) examined the nonlinear stability of triangular
points. Moser gave some modifications in Arnold's theorem. Then
Deprit and Deprit (1967) investigated the nonlinear stability of
triangular points by applying Moser's modified version of Arnold's
theorem (1961).

Bhatnagar and Hallan (1983) studied the effect of perturbations on
the nonlinear stability of triangular points. Maciejewski and
Gozdziewski (1991) described the normalization algorithms of
Hamiltonian near an equilibrium point. Niedzielska (1994) investigated the nonlinear
stability of the libration points in the photogravitational
restricted three body problem. Mishra P. and Ishwar B.(1995) studied second order normalization in the generalized restricted problem of three bodies, smaller primary being an oblate spheroid. Ishwar B.(1997) studied nonlinear stability in the generalized restricted three body problem.

In this paper we have performed Birkhoff's normalization of the Hamiltonian. For this we have utilised Henrard's  method and expanded the coordinates of the third body in Double d'Alembert series. We have found the values of first and second order components. The second order components are obtained as solutions of the two partial differential equations. We have employed the first condition of KAM theorem in solving these equations. The first and second order components are affected by radiation pressure, oblateness and P-R drag. Finaly we obtained the third order part $H_3$ of the Hamiltonian in $I_1^{1/2}I_2^{1/2}$ zero.

\section{Location of Triangular Equilibrium Points}
\newcommand{\abc}{{\Bigl(1+\frac{5}{2}A_2\Bigr) }}
\newcommand{\aac}{{(1-\frac{A_2}{2}) }}
\newcommand{\amc}{{(1-\mu) }}
\newcommand{\adc}{{\frac{\delta ^2}{2} }}
\newcommand{\zab}{{\Bigl(1+\frac{5A_2}{2r^2_2}\Bigr)}}
\newcommand{\zw}{{\frac{W_1}{r^2_1}}}
\newcommand{\zwf}{{\frac{W_1}{r^4_1}}}
\newcommand{\zx}{{(x+\mu)}}
\newcommand{\zox}{{(x+\mu-1)}}
\newcommand{\zd}{{\displaystyle}}
\newcommand{\zabs}{{\Bigl(1+\frac{5A_2}{2{r^2_2}_*}\Bigr)}}
\newcommand{\zws}{{\frac{W_1}{{r^2_1}_*}}}
\newcommand{\zwfs}{{\frac{W_1}{{r^4_1}_*}}}
\newcommand{\zxs}{{(x_*+\mu)}}
\newcommand{\zoxs}{{(x_*+\mu-1)}}
\newcommand{\az}{{1,0}}
\newcommand{\za}{{0,1}}
\newcommand{\zae}{A_2\epsilon}
\newcommand{\zwe}{nW_1\epsilon}
Equations of motions are
\begin{align}
\ddot{x}-2n\dot{y}&=U_x ,\quad\text{where},\quad U_x=\frac{\partial{U_1}}{\partial{x}}-\frac{W_{1}N_1}{r^2_1}\\
\ddot{y}+2n\dot{x}&=U_y,\hspace{.85in}U_y=\frac{\partial{U_1}}{\partial{y}}-\frac{W_{1}N_2}{r^2_1}\\
U_1&=\zd{\frac{n^2(x^2+y^2)}{2}}+\frac{\amc{q_1}}{r_1}+\frac{\mu}{r_2}+\frac{\mu{A_2}}{2r^3_2}
\end{align}
\begin{gather*}
r^2_1=\zx^2+y^2,\quad  r^2_2=\zox^2+y^2,\quad n^2=1+\frac{3}{2}A_2,\\
N_1=\frac{\zx[\zx\dot{x}+y\dot{y}]}{r^2_1}+\dot{x}-ny,\quad
N_2=\frac{y[\zx\dot{x}+y\dot{y}]}{r^2_1}+\dot{y}+n\zx
\end{gather*}
$W_1=\frac{(1-\mu)(1-q_1)}{c_d}$,
$\mu=\frac{m_2}{m_1+m_2}\leq\frac{1}{2}$, $m_1,m_2$ be  the  masses of
the primaries, $A_2=\frac{r^2_e-r^2_p}{5r^2}$ be the oblateness
coefficient, $r_e$ and$r_p$ be the equatorial and polar radii
respectively $r$ be the distance between primaries,
$q=\bigl(1-\frac{F_p}{F_g}\bigr)$ be the mass reduction factor
expressed in terms of the particle's radius $a$, density $\rho$ and
radiation pressure efficiency factor $\chi$ (in the C.G.S.system)
i.e., $q=1-\zd{\frac{5.6\times{10^{-5}}\chi}{a\rho}}$. Assumption
$q=constant$ is equivalent to neglecting fluctuation in the beam of
solar radiation and the effect of solar radiation, the effect of the
planet's shadow, obviously $q\leq1$. Triangular equilibrium points
are given by $U_x=0,U_y=0,z=0,y\neq{0}$, then we have
\begin{align}
 x_*&=x_0\Biggl\{1-\zd{\frac{nW_1\bigl[\amc\abc+\mu\aac\adc\bigr]}{3\mu\amc{y_0 x_0}}}-\adc\frac{A_2}{x_0}\Biggr\} \label{eq:1x}\\
y_*&=y_0{\Biggl\{1-\zd{\frac{nW_1\delta^2\bigl[2\mu-1-\mu(1-\frac{3A_2}{2})\adc+7\amc\frac{A_2}{2}\bigr]}{3\mu\amc{y^3_0}}}-\zd{\frac{\delta^2\bigl(1-\adc)A_2}{y^2_0}}\Biggr\}^{1/2}
}
\end{align}
where
$x_0=\adc-\mu$, $y_0=\pm\delta\bigl(1-\frac{\delta^2}{4}\bigr)^{1/2}$
and $\delta=q^{1/3}_1$, as in preprint Kushvah \& Ishwar (2006)
\section{Second Order Normalization}
The Lagrangian function of the problem can be written as
\begin{align}
L&=\frac{1}{2}(\dot{x}^2+\dot{y}^2)+n(x\dot{y}-\dot{x}y)+\frac{n^2}{2}(x^2+y^2)+\frac{\amc{q_1}}{r_1}+\frac{\mu}{r_2}+\frac{\mu{A_2}}{2r^3_2}\\\notag
&+W_1\Bigl\{\frac{\zx\dot{x}+y\dot{y}}{2r^2_1}-n
\arctan{\frac{y}{\zx}}\Bigr\}\\\notag
\end{align}
and the Hamiltonian  is $H=-L+P_x\dot{x}+P_y\dot{y}$, where
$P_x,P_y$ are the momenta coordinates given by \[
P_x=\frac{\partial{L}}{\partial{\dot{x}}}=\dot{x}-ny+\frac{W_1}{2r_1^2}\zx,
\quad
P_y=\frac{\partial{L}}{\partial{\dot{y}}}=\dot{y}+nx+\frac{W_1}{2r_1^2}y
\]
For simplicity we suppose  $q_1=1-\epsilon$, with $|\epsilon|<<1$
then coordinates of triangular equilibrium points  can be
written in the form
\begin{align}
x&=\frac{\gamma}{2}-\frac{\epsilon}{3}-\frac{A_2}{2}+\frac{A_2
\epsilon}{3}-\frac{(9+\gamma)}{6\sqrt{3}}nW_1-\frac{4\gamma
\epsilon}{27\sqrt{3}}nW_1 \\
y&=\frac{\sqrt{3}}{2}\Bigl\{1-\frac{2\epsilon}{9}-\frac{A_2}{3}-\frac{2A_2
\epsilon}{9}+\frac{(1+\gamma)}{9\sqrt{3}}nW_1-\frac{4\gamma
\epsilon}{27\sqrt{3}}nW_1\Bigr\}
\end{align}
where $\gamma=1-2\mu$.
 We shift the origin to $L_4$. For that, we change
$x\rightarrow {x_*}+x$ and  $y\rightarrow{y_*}+y$. Let $a=x_*+\mu, b=y_*$ so
that
\begin{align}
a&= \frac{1}{2} \Biggl\{-\frac{2\epsilon}{3}-A_2+\frac{2A_2
\epsilon}{3}-\frac{(9+\gamma)}{3\sqrt{3}}nW_1-\frac{8\gamma
\epsilon}{27\sqrt{3}}nW_1 \bigr\}\\
b&=\frac{\sqrt{3}}{2}\Bigl\{1-\frac{2\epsilon}{9}-\frac{A_2}{3}-\frac{2A_2
\epsilon}{9}+\frac{(1+\gamma)}{9\sqrt{3}}nW_1-\frac{4\gamma
\epsilon}{27\sqrt{3}}nW_1\Bigr\}
\end{align}
Expanding $L$ in power series of $x\mbox{and}y$, we get
\begin{eqnarray}
 L&=&L_0+L_1+L_2+L_3+\cdots \\
H&=&H_0+H_1+H_2+H_3+\cdots =-L+P_x{\dot{x}}+P_y{\dot{y}}
  \end{eqnarray}
  where $L_0,L_1,L_2,L_3 \ldots$ are constants, first order term,
  second order term, \ldots respectively.
  Third order term $H_3$ of Hamiltonian can be written as
  \begin{equation}\label{eq:h3}
  H_3=-L_3=-\frac{1}{3!}\left\{x^3T_1+3x^2yT_2+3xy^2T_3+y^3T_4+6T_5\right\}
  \end{equation}
  where
 \begin{eqnarray} T_1&=&\frac{3}{16}\biggl[\frac{16}{3}\epsilon+6A_2-\frac{979}{18}\zae+\frac{(143+9\gamma)}{6\sqrt{3}}nW_1+\frac{(459+376\gamma)}{27\sqrt{3}}\zwe\\&&+\gamma\left\{14+\frac{4\epsilon}{3}+25A_2-\frac{1507 }{18}\zae-\frac{(215+29\gamma)}{6\sqrt{3}}nW_1
-\frac{2(1174+169\gamma)}{27\sqrt{3}}\zwe\right\}\biggr]\notag\\
T_2&=&\frac{3\sqrt{3}}{16}\biggl[14-\frac{16}{3}\epsilon+\frac{A_2}{3}-\frac{367}{18}\zae+\frac{115(1+\gamma)}{18\sqrt{3}}nW_1-\frac{(959-136\gamma)}{27\sqrt{3}}\zwe\\&&+\gamma\left\{\frac{32\epsilon}{3}+40A_2-\frac{382}{9}\zae+\frac{(511+53\gamma)}{6\sqrt{3}}nW_1-\frac{(2519-24\gamma)}{27\sqrt{3}}\zwe\right\}\biggr]\notag\\
 T_3&=&\frac{-9}{16}\biggl[\frac{8}{3}\epsilon+\frac{203A_2}{6}-\frac{625}{54}\zae-\frac{(105+15\gamma)}{18\sqrt{3}}nW_1-\frac{(403-114\gamma)}{81\sqrt{3}}\zwe\\&&+\gamma\left\{2-\frac{4\epsilon}{9}+\frac{55A_2}{2}-\frac{797}{54}\zae+\frac{(197+23\gamma)}{18\sqrt{3}}nW_1-\frac{(211-32\gamma)}{81\sqrt{3}}\zwe\right\}\biggr]\notag\\
 T_4&=&\frac{-9\sqrt{3}}{16}\biggl[2-\frac{8}{3}\epsilon+\frac{23A_2}{3}-44\zae-\frac{(37+\gamma)}{18\sqrt{3}}nW_1-\frac{(219+253\gamma)}{81\sqrt{3}}\zwe\\&&+\gamma\left\{4\epsilon+\frac{88}{27}\zae+\frac{(241+45\gamma)}{18\sqrt{3}}nW_1-\frac{(1558-126\gamma)}{81\sqrt{3}}\zwe\right\}\biggr]\notag\\ T_5&=&\frac{W_1}{2(a^2+b^2)^3}\biggl[(a\dot{x}+b\dot{y})\left\{3(ax+by)-(bx-ay)^2\right\}-2(x\dot{x}+y\dot{y})(ax+by)(a^2+b^2)\biggr]\label{eq:t5}\end{eqnarray}

  In order to perform Birkhoff's normalization, we use Henrard's
  method[Deprit and Deprit Brtholom\'{e} (1967)] for which the
  coordinates $(x,y)$ of infinitesimal body, to be expanded in
  double d'Alembert series  $x=\sum_{n\geq1}B_n^\az,\quad y=\sum_{n\geq 1}B_n^\za$
  where the homogeneous components $B_n^\az$ and $B_n^\za$ of
  degree $n$ are of the form
  \begin{equation}
\sum_{0\leq{m}\leq{n}}
I_1^{\frac{n-m}{2}}I_2^{\frac{m}{2}}\sum_{(p,q)}C_{n-m,m,p,q}
\cos{(p\phi_1+q\phi_2)}+S_{n-m,m,p,q} \sin{(p\phi_1+q\phi_2)}
  \end{equation}
  The condition in double summation are (i) $p$ runs over those
  integers in the interval $0\leq p\leq n-m$ that have the same
  parity as $n-m$ (ii) $q$ runs over those integers in the interval $-m\leq q\leq
  m$ that have the same parity as $m$. Here $I_1$, $I_2$ are the
  action momenta coordinates which are to be taken as constants of
  integer, $\phi_1$, $\phi_2$ are angle coordinates to be
  determined as linear functions of time in such a way that $\dot\phi_1=\omega_1+\sum_{n\geq 1}f_{2n}(I_1,I_2),\dot\phi_2=\omega_2+\sum_{n\geq 1}g_{2n}(I_1,I_2)$  where  $\omega_1,\omega_2$ are the basic  frequencies, $f_{2n}$ and  $g_{2n}$ are of the form
  \begin{eqnarray}
  f_{2n}&=&\sum_{0\leq m\leq n}{f'}_{2(n-m),2m}I_1^{n-m}I_2^m\\
  g_{2n}&=&\sum_{0\leq m\leq n}{g'}_{2(n-m),2m}I_1^{n-m}I_2^m
  \end{eqnarray}
The first order components $B_1^\az$ and $B_1^\za$  in
$I_1$, $I_2$ are the values of $x$ and $y$ given by
\begin{equation}  X=JT \quad  \text{where}  \quad X=\left[\begin{array}{c}
x\\y\\p_x\\p_y\end{array}\right],J=[J_{ij}]_{1\leq i\leq j \leq 4},\
T=\left[\begin{array}{c}
Q_1\\Q_2\\p_1\\p_2\end{array}\right]
\end{equation}
\begin{equation}
P_i= (2 I_i\omega_i)^{1/2}\cos{\phi_i},  \quad 
Q_i= (\frac{2 I_i}{\omega_i})^{1/2}\sin{\phi_i}, \quad (i=1,2)
\end{equation}
\begin{eqnarray}\label{eq:b1az} B_1^\az&=&J_{13}\sqrt{2\omega_1I_1}\cos{\phi_1}+J_{14}\sqrt{2\omega_2I_2}\cos{\phi_2}\\
B_1^\za&=&J_{21}\sqrt{\frac{2I_1}{\omega_1}}\sin{\phi_1}+J_{22}\sqrt{\frac{2I_2}{\omega_2}}\sin{\phi_2}+J_{23}\sqrt{2I_1}{\omega_1}\cos{\phi_1}+J_{24}\sqrt{2I_2}{\omega_2}\sin{\phi_2}
\end{eqnarray}
where \begin{eqnarray}
J_{13}&=&\frac{l_1}{2\omega_1k_1}\left\{1-\frac{1}{2l_1^2}\left[\epsilon+\frac{45A_2}{2}-\frac{717A_2\epsilon}{36}+\frac{(67+19\gamma)}{12\sqrt{3}}nW_1
-\frac{(431-3\gamma)}{27\sqrt{3}}nW_1\epsilon\right]\right.\\&&+\frac{\gamma}{2l_1^2}\left[3\epsilon-\frac{29A_2}{36}-\frac{(187+27\gamma)}{12\sqrt{3}}nW_1
-\frac{2(247+3\gamma)}{27\sqrt{3}}nW_1\epsilon\right]\notag\\&&-\frac{1}{2k_1^2}\left[\frac{\epsilon}{2}-3A_2-\frac{73A_2\epsilon}{24}+\frac{(1-9\gamma)}{24\sqrt{3}}nW_1
+\frac{(53-39\gamma)}{54\sqrt{3}}nW_1\epsilon\right]\notag\\&&-\frac{\gamma}{4k_1^2}\left[\epsilon-3A_2-\frac{299A_2\epsilon}{72}-\frac{(6-5\gamma)}{12\sqrt{3}}nW_1
-\frac{(266-93\gamma)}{54\sqrt{3}}nW_1\epsilon\right]\notag\\&&\left.+\frac{\epsilon}{4l_1^2k_1^2}\left[\frac{3A_2}{4}
+\frac{(33+14\gamma)}{12\sqrt{3}}nW_1\right]+\frac{\gamma\epsilon}{8l_1^2k_1^2}\left[\frac{347A_2}{36}
-\frac{(43-8\gamma)}{4\sqrt{3}}nW_1 \right]\right\}\notag\end{eqnarray}
 \begin{eqnarray} J_{14}&=&\frac{l_2}{2\omega_2k_2}\left\{1-\frac{1}{2l_2^2}\left[\epsilon+\frac{45A_2}{2}-\frac{717A_2\epsilon}{36}+\frac{(67+19\gamma)}{12\sqrt{3}}nW_1
-\frac{(431-3\gamma)}{27\sqrt{3}}nW_1\epsilon\right]\right.\\&&-\frac{\gamma}{2l_2^2}\left[3\epsilon-\frac{293A_2}{36}+\frac{(187+27\gamma)}{12\sqrt{3}}nW_1
-\frac{2(247+3\gamma)}{27\sqrt{3}}nW_1\epsilon\right]\notag\\&&-\frac{1}{2k_2^2}\left[\frac{\epsilon}{2}-3A_2-\frac{73A_2\epsilon}{24}+\frac{(1-9\gamma)}{24\sqrt{3}}nW_1
+\frac{(53-39\gamma)}{54\sqrt{3}}nW_1\epsilon\right]\notag\\&&+\frac{\gamma}{2k_2^2}\left[\epsilon-3A_2-\frac{299A_2\epsilon}{72}-\frac{(6-5\gamma)}{12\sqrt{3}}nW_1
-\frac{(268-9\gamma)}{54\sqrt{3}}nW_1\epsilon\right]\notag\\&&\left.-\frac{\epsilon}{4l_2^2k_2^2}\left[\frac{33A_2}{4}
+\frac{(1643-93\gamma)}{216\sqrt{3}}nW_1\right]+\frac{\gamma\epsilon}{4l_2^2k_2^2}\left[\frac{737A_2}{72}
-\frac{(13+2\gamma)}{\sqrt{3}}nW_1 \right]\right\}\notag\end{eqnarray}
\begin{eqnarray} J_{21}&=&-\frac{4n\omega_1}{l_1k_1}\left\{1+\frac{1}{2l_1^2}\left[\epsilon+\frac{45A_2}{2}-\frac{717A_2\epsilon}{36}+\frac{(67+19\gamma)}{12\sqrt{3}}nW_1
-\frac{(413-3\gamma)}{27\sqrt{3}}nW_1\epsilon\right]\right.\\&&-\frac{\gamma}{2l_1^2}\left[3\epsilon-\frac{293A_2}{36}+\frac{(187+27\gamma)}{12\sqrt{3}}nW_1
-\frac{2(247+3\gamma)}{27\sqrt{3}}nW_1\epsilon\right]\notag\\&&-\frac{1}{2k_1^2}\left[\frac{\epsilon}{2}-3A_2-\frac{73A_2\epsilon}{24}+\frac{(1-9\gamma)}{24\sqrt{3}}nW_1
+\frac{(53-39\gamma)}{54\sqrt{3}}nW_1\epsilon\right]\notag\\&&-\frac{\gamma}{4k_1^2}\left[\epsilon-3A_2-\frac{299A_2\epsilon}{72}-\frac{(6-5\gamma)}{12\sqrt{3}}nW_1
-\frac{(268-93\gamma)}{54\sqrt{3}}nW_1\epsilon\right]\notag\\&&\left.+\frac{\epsilon}{8l_1^2k_1^2}\left[\frac{33A_2}{4}+\frac{(68-10\gamma)}{24\sqrt{3}}nW_1\right]+\frac{\gamma\epsilon}{8l_1^2k_1^2}\left[\frac{242A_2}{9}
+\frac{(43-8\gamma)}{4\sqrt{3}}nW_1 \right]\right\}\notag\end{eqnarray}
\begin{eqnarray} J_{22}&=&\frac{4n\omega_2}{l_2k_2}\left\{1+\frac{1}{2l_2^2}\left[\epsilon+\frac{45A_2}{2}-\frac{717A_2\epsilon}{36}+\frac{(67+19\gamma)}{12\sqrt{3}}nW_1
-\frac{(413-3\gamma)}{27\sqrt{3}}nW_1\epsilon\right]\right.\\&&-\frac{\gamma}{2l_2^2}\left[3\epsilon-\frac{293A_2}{36}+\frac{(187+27\gamma)}{12\sqrt{3}}nW_1
-\frac{2(247+3\gamma)}{27\sqrt{3}}nW_1\epsilon\right]\notag\\&&+\frac{1}{2k_2^2}\left[\frac{\epsilon}{2}-3A_2-\frac{73A_2\epsilon}{24}+\frac{(1-9\gamma)}{24\sqrt{3}}nW_1
+\frac{(53-39\gamma)}{54\sqrt{3}}nW_1\epsilon\right]\notag\\&&-\frac{\gamma}{4k_2^2}\left[\epsilon-3A_2-\frac{299A_2\epsilon}{72}-\frac{(6-5\gamma)}{12\sqrt{3}}nW_1
-\frac{(268-93\gamma)}{54\sqrt{3}}nW_1\epsilon\right]\notag\\&&\left.+\frac{\epsilon}{4l_2^2k_2^2}\left[\frac{33A_2}{4}+\frac{(34+5\gamma)}{12\sqrt{3}}nW_1
\right]+\frac{\gamma\epsilon}{8l_2^2k_2^2}\left[\frac{75A_2}{2}+\frac{(43-8\gamma)}{4\sqrt{3}}nW_1
 \right]\right\}\notag\end{eqnarray}
\begin{eqnarray} J_{23}&=&\frac{\sqrt{3}}{4\omega_1l_1k_1}\left\{2\epsilon+6A_2+\frac{37A_2\epsilon}{2}-\frac{(13+\gamma)}{2\sqrt{3}}nW_1
+\frac{2(79-7\gamma)}{9\sqrt{3}}nW_1\epsilon\right.\\&&-\gamma\left[6+\frac{2\epsilon}{3}+13A_2-\frac{33A_2\epsilon}{2}+\frac{(11-\gamma)}{2\sqrt{3}}nW_1
-\frac{(186-\gamma)}{9\sqrt{3}}nW_1\epsilon\right]\notag\\&&+\frac{1}{2l_1^2}\left[51A_2+\frac{(14+8\gamma)}{3\sqrt{3}}nW_1\right]-\frac{\epsilon}{k_1^2}\left[3A_2
+\frac{(19+6\gamma)}{6\sqrt{3}}nW_1\right]\notag\\&&-\frac{\gamma}{2l_1^2}\left[6\epsilon+135A_2-\frac{808A_2\epsilon}{9}-\frac{(67+19\gamma)}{2\sqrt{3}}nW_1
-\frac{(755+19\gamma)}{9\sqrt{3}}nW_1\epsilon\right]\notag\\&&-\frac{\gamma}{2k_1^2}\left[3\epsilon-18A_2-\frac{55A_2\epsilon}{4}-\frac{(1-9\gamma)}{4\sqrt{3}}nW_1
+\frac{(923-60\gamma)}{12\sqrt{3}}nW_1\epsilon\right]\notag\\&&\left.+\frac{\gamma\epsilon}{8l_1^2k_1^2}\left[\frac{9A_2}{2}
+\frac{(34-5\gamma)}{2\sqrt{3}}nW_1\right]\right\}\qquad\notag\end{eqnarray}
\begin{eqnarray}				 J_{24}&=&\frac{\sqrt{3}}{4\omega_2l_2k_2}\left\{2\epsilon+6A_2+\frac{37A_2\epsilon}{2}-\frac{(13+\gamma)}{2\sqrt{3}}nW_1
+\frac{2(79-7\gamma)}{9\sqrt{3}}nW_1\epsilon\right.\label{eq:j24}\\ &&-\gamma\left[6+\frac{2\epsilon}{3}+13A_2-\frac{33A_2\epsilon}{2}+\frac{(11-\gamma)}{2\sqrt{3}}nW_1
-\frac{(186-\gamma)}{9\sqrt{3}}nW_1\epsilon\right]\notag\\&&-\frac{1}{2l_2^2}\left[51A_2+\frac{(14+8\gamma)}{3\sqrt{3}}nW_1\right]-\frac{\epsilon}{k_2^2}\left[3A_2
+\frac{(19+6\gamma)}{6\sqrt{3}}nW_1\right]\notag\\&&-\frac{\gamma}{2l_2^2}\left[6\epsilon+135A_2-\frac{808A_2\epsilon}{9}-\frac{(67+19\gamma)}{2\sqrt{3}}nW_1
-\frac{(755+19\gamma)}{9\sqrt{3}}nW_1\epsilon\right]\notag\\&&-\frac{\gamma}{2k_1^2}\left[3\epsilon-18A_2-\frac{55A_2\epsilon}{4}-\frac{(1-9\gamma)}{4\sqrt{3}}nW_1
+\frac{(923-60\gamma)}{12\sqrt{3}}nW_1\epsilon\right]\notag\\&&\left.-\frac{\gamma\epsilon}{4l_1^2k_1^2}\left[\frac{99A_2}{2}
+\frac{(34-5\gamma)}{2\sqrt{3}}nW_1\right]\right\}\notag\end{eqnarray}
with $l_j^2=4\omega_j^2+9,(j=1,2)$ and $ k_1^2=2\omega_1^2-1, k_2^2=-2\omega_2^2+1 $. 
In order to findout  the second order components $B_2^\az,B_2^\za$ we
consider  Lagrange's  equations of motion
\begin{equation}
\frac{d}{dt}(\frac{\partial L}{\partial \dot x })-\frac{\partial L}{\partial x }=0, \quad
\frac{d}{dt}(\frac{\partial L}{\partial \dot y })-\frac{\partial
L}{\partial y }=0 \end{equation}
\begin{equation}
\text{i.e.}\quad \left.\begin{array}{l c l}
\ddot x-2n\dot y+(2E-n^2)x+Gy&=&\frac{\partial L_3}{\partial x }+\frac{\partial L_4}{\partial x }\\
&&\\\ddot x+2n\dot x+(2F-n^2)y+Gx&=&\frac{\partial L_3}{\partial y
}+\frac{\partial L_4}{\partial y }\end{array}
\right\}\label{eq:lgeq}\end{equation}
Since $x$ and $y$ are double
d'Alembert series, $x^jx^k(j\geq0,k\geq0,j+k\geq0)$ is also a double
d'Alembert series, the time derivatives $\dot x ,\dot y ,\ddot x,
\ddot y $ are also double d'Alembert series. We can write
\[\dot x=\sum_{n\geq 1} \dot x_n, \quad\dot y=\sum_{n\geq 1} \dot
y_n,\quad\ddot x=\sum_{n\geq 1} \ddot x_n,\quad \ddot y=\sum_{n\geq
1} \ddot y_n \] where $\dot x ,\dot y ,\ddot x, \ddot y $ are
homogeneous components of degree $n$ in $I_1^{1/2},I_2^{1/2}$ i.e.
\begin{eqnarray} \dot x &=&
\frac{d}{dt}\sum_{n\geq 1}B_n^\az=\sum_{n\geq 1}\left[\frac{\partial
B_n^\az}{\partial{\phi_1}}(\omega_1+f_2+f_4+\cdots)+\frac{\partial
B_n^\az}{\partial{\phi_2}}(-\omega_2+g_2+g_4+\cdots)\right]\end{eqnarray}
We write three components $\dot x_1 ,\dot x_2 ,\dot x_3$ of $\dot x$
\begin{eqnarray}
\dot x_1&=&\omega_1\frac{\partial
B_1^\az}{\partial{\phi_1}}-\omega_2\frac{\partial
B_1^\az}{\partial{\phi_2}}=DB_1^\az\\
 \dot
x_2&=&\omega_1\frac{\partial
B_2^\az}{\partial{\phi_1}}-\omega_2\frac{\partial
B_2^\az}{\partial{\phi_2}}=DB_2^\az\\
\dot x_3&=&\omega_1\frac{\partial B_3^\az}{\partial{\phi_1}}-\omega_2\frac{\partial
B_3^\az}{\partial{\phi_2}}+f_2\frac{\partial
B_1^\az}{\partial{\phi_1}}-g_2\frac{\partial
B_1^\az}{\partial{\phi_2}}\notag\\
&=&DB_2^\az+f_2\frac{\partial
B_1^\az}{\partial{\phi_1}}-g_2\frac{\partial
B_1^\az}{\partial{\phi_2}}
\end{eqnarray}
where \begin{equation}D\equiv \omega_1\frac{\partial\
}{\partial{\phi_1}}-\omega_2\frac{\partial\
}{\partial{\phi_2}}\end{equation}
 Similarly three components $\ddot
x_1 ,\ddot x_2 ,\ddot x_3$ of $\ddot x$ are
\begin{eqnarray*}
\ddot x_1 &=&D^2B_1^\az, \quad \ddot x_2=D^2B_2^\az,\quad \ddot
x_3=D^2B_3^\az+2\omega_1f_2\frac{\partial^2B_1^\az}{\partial\phi_1^2}-2\omega_2g_2\frac{\partial^2B_1^\az}{\partial\phi_2^2}
\end{eqnarray*}
In similar manner we can write the components of $\dot y, \ddot y$.
Putting the values of  $x, y, \dot x ,\dot y ,\ddot x $ and  $\ddot
y$ in terms of double d'Alembert series in equation
(~\ref{eq:lgeq}) we get
\begin{equation}
\left(D^2+2E-1-\frac{3}{2}A_2\right)B_2^\az-\left\{2\left(
1+\frac{3}{4}A_2\right)D-G\right\}B_2^\za=X_2 \label{eq:x2}
\end{equation}
\begin{equation}\label{eq:y2}
\left\{2\left(
1+\frac{3}{4}A_2\right)D+G\right\}B_2^\az+\left(D^2+2F-1-\frac{3}{2}A_2\right)B_2^\za=Y_2
\end{equation} where \[X_2=\left[\frac{\partial
L_3}{\partial x}\right]_{x=B_1^\az,y=B_1^\za} \quad \text{and} \quad
Y_2=\left[\frac{\partial L_3}{\partial
y}\right]_{x=B_1^\az,y=B_1^\za}\] 
These are two simultaneous partial differential equations in $B_2^\az$ and $B_2^\za$. We solve these
equations to find the values of $B_2^\az$ and $B_2^\za$, from Eq.
(~\ref{eq:x2}) and (~\ref{eq:y2})
\begin{equation}
\triangle_1 \triangle_1B_2^\az=\Phi_2, \quad \triangle_1 \triangle_1B_2^\za=-\Psi_2 \label{eq:phi_si} \quad
\text{where} \quad \triangle_1=D^2+\omega_1^2, \triangle_2=D^2+\omega_2^2 \end{equation}
\begin{equation}
\Phi_2=(D^2+2F-n^2)X_2+(2nD-G)Y_2 \label{eq:phi2}
\end{equation}
\begin{equation}
\Psi_2=(2nD+G)X_2-(D^2+2E-n^2)Y_2 \label{eq:psi2}
\end{equation}
The Eq.(~\ref{eq:phi_si}) can  be solved for $B_2^\az$
and $B_2^\za$ by putting the formula
\[\frac{1}{\triangle_1\triangle_2}\left\{\begin{array}{c}\cos(p\phi_1+q\phi_2)\\ \mbox{or} \\\sin(p\phi_1+q\phi_2)\end{array}=\frac{1}{\triangle_1\triangle_2}\left\{\begin{array}{c}\cos(p\phi_1+q\phi_2)\\\mbox{or} \\\sin(p\phi_1+q\phi_2)\end{array}\right.\right.\]
where \[\triangle_{p,q}=\left[
\omega_1^2-(\omega_1p-\omega_2q)^2\right]\left[
\omega_2^2-(\omega_1p-\omega_2q)^2\right]
\] 
provided $\triangle_{p,q}\neq0$.
Since $\triangle_{1,0}=0, \triangle_{0,1}=0$ the terms
$\cos\phi_1,\sin\phi_1,\cos\phi_2,\sin\phi_2$ are the critical terms.  
$\phi$ and $\psi$ are free from such terms. By condition(1) of
Moser's theorem $k_1\omega_1+k_2\omega_2\neq 0$  for all pairs
$(k_1,k_2)$ of integers such that $|k_1|+|k_2|\leq4$, therefore each
of $\omega_1, \omega_2, \omega_1\pm2\omega_2,\omega_2\pm2\omega_1$
is different from zero and consequently none of the divisors
$\triangle_{0,0}, \triangle_{0,2}, \triangle_{2,0}, \triangle_{1,1}, \triangle_{1,-1}$ is zero. The second
order components $B_2^\az, B_2^\za$ are as follows:
\begin{eqnarray}\label{eq:b2az}
B_2^\az&=&r_1I_1+r_2I_2+r_3I_1\cos2\phi_1+r_4I_2\cos2\phi_2+r_5I_1^{1/2}I_2^{1/2}\cos(\phi_1-\phi_2)\\&&
+r_6I_1^{1/2}I_2^{1/2}\cos(\phi_1+\phi_2)+r_7I_1\sin2\phi_1+r_8I_2\sin2\phi_2\notag\\
&&+r_9I_1^{1/2}I_2^{1/2}\sin(\phi_1-\phi_2)+r_{10}I_1^{1/2}I_2^{1/2}\sin(\phi_1+\phi_2)\notag \end{eqnarray}
and\begin{eqnarray}
B_2^\za&=&-\left\{s_1I_1+s_2I_2+s_3I_1\cos2\phi_1+s_4I_2\cos2\phi_2+s_5I_1^{1/2}I_2^{1/2}\cos(\phi_1-\phi_2)\right.\\
&&
+s_6I_1^{1/2}I_2^{1/2}\cos(\phi_1+\phi_2)+s_7I_1\sin2\phi_1+s_8I_2\sin2\phi_2\notag\\
&&+\left.s_9I_1^{1/2}I_2^{1/2}\sin(\phi_1-\phi_2)+s_{10}I_1^{1/2}I_2^{1/2}\sin(\phi_1+\phi_2)\right\}\notag\label{eq:b2za}\notag
\end{eqnarray}
where

\begin{eqnarray}
&r_1=&\frac{1}{\omega_1^2\omega_2^2}\left\{
J_{13}^2\omega_1F_4+J_{13}J_{23}\omega_1F_4'+\left(\frac{J_{21}^2}{\omega_1}+J_{23}^2\omega_1\right)F_4''\right\}
\end{eqnarray}
\begin{eqnarray}
&r_2=&\frac{1}{\omega_1^2\omega_2^2}\left\{
J_{14}^2\omega_2F_4+J_{14}J_{24}\omega_2F_4'+\left(\frac{J_{22}^2}{\omega_2}+J_{24}^2\omega_2\right)F_4{''}\right\}
\end{eqnarray}
\begin{eqnarray}
&r_3=&\frac{-1}{3\omega_1^2(4\omega_1^2-\omega_2^2)}\Biggl\{
8\omega_1^3J_{21}(J_{13}F_1'+2J_{23}F_1'')+4\omega_1^2\biggl[(J_{13}F_2+J_{23}F_2'')J_{13}\omega_1-\biggl(\frac{J_{21}^2}{\omega_1}-J_{23}^2\omega_1\biggr)F_1''\biggr]\\
&&-2\omega_1J_{21}(J_{13}F_3'+2J_{23}F_3'')-\omega_1J_{13}(J_{13}F_4+J_{23}F_4'')\omega_1
+\biggl(\frac{J_{21}^2}{\omega_1}-J_{23}^2\omega_1\biggr)F_1''\Biggr\}\notag\\
&&\notag
\end{eqnarray}
\begin{eqnarray}
&r_4=&\frac{1}{3\omega_2^2(4\omega_2^2-\omega_1^2)}\Biggl\{8\omega_2^3J_{22}(J_{14}F_1'+2J_{24}F_1'')-4\omega_2^2\biggl[(J_{14}F_2+J_{24}F_2'')J_{14}\omega_2-\biggl(\frac{J_{22}^2}{\omega_2}-J_{24}^2\omega_2\biggr)F_2''\biggr]\\&&-2\omega_2J_{22}(J_{14}F_3'+2J_{24}F_3'')-\omega_2J_{14}(J_{14}F_4+J_{24}F_4'')\omega_2-\biggl(\frac{J_{22}^2}{\omega_2}-J_{24}^2\omega_2\biggr)F_4''\Biggr\}\notag\end{eqnarray}
\begin{eqnarray}
&r_5&=\frac{1}{\omega_1\omega_2(2\omega_1+\omega_2)(4\omega_1+2\omega_2)}\Biggl\{(\omega_1+\omega_2)^3\biggl[\bigl\{J_{13}J_{22}(\frac{\omega_1}{\omega_2})^{1/2}-J_{14}J_{21}(\frac{\omega_2}{\omega_1})^{1/2}\bigr\}F_1'\\&&-2\bigl\{J_{21}J_{24}(\frac{\omega_2}{\omega_1})^{1/2}-J_{22}J_{23}(\frac{\omega_1}{\omega_2})^{1/2}\bigr\}F_1''\biggl]-(\omega_1+\omega_2)^2\biggl[\bigr\{2\bigl\{J_{13}J_{14}F_2\notag\\&&+(J_{13}J_{24}+J_{14}J_{23})F_2'\bigr\}(\omega_1\omega_2)^{1/2}+\bigr\{\frac{J_{21}J_{22}}{(\omega_1\omega_2)^{1/2}}+J_{23}J_{24}(\omega_1\omega_2)^{1/2}\bigr\}F_2''\biggr]\notag
\\
&&-(\omega_1+\omega_2)\biggl[\bigl\{J_{13}J_{22}(\frac{\omega_1}{\omega_2})^{1/2}-J_{14}J_{21}(\frac{\omega_2}{\omega_1})^{1/2}\bigr\}F_3'-2\bigl\{J_{21}J_{24}(\frac{\omega_2}{\omega_1})^{1/2}-J_{22}J_{23}(\frac{\omega_1}{\omega_2})^{1/2}\bigr\}F_3''\biggl]\notag\\&&+\biggl[\bigr\{2\bigl\{J_{13}J_{14}F_4+(J_{13}J_{24}+J_{14}J_{23})F_4'\bigr\}(\omega_1\omega_2)^{1/2}+2\bigr\{\frac{J_{21}J_{22}}{(\omega_1\omega_2)^{1/2}}+J_{23}J_{24}(\omega_1\omega_2)^{1/2}\bigr\}F_4''\biggr]\Biggr\}\notag
\end{eqnarray}
\begin{eqnarray}
&r_6=&\frac{-1}{\omega_1\omega_2(2\omega_1-\omega_2)(4\omega_1-2\omega_2)}
\Biggl\{(\omega_1-\omega_2)^3\biggl[\bigl\{J_{13}J_{22}(\frac{\omega_1}{\omega_2})^{1/2}-J_{14}J_{21}(\frac{\omega_2}{\omega_1})^{1/2}\bigr\}F_1'\\
&&+2\bigl\{J_{21}J_{24}(\frac{\omega_2}{\omega_1})^{1/2}+J_{22}J_{23}(\frac{\omega_1}{\omega_2})^{1/2}\bigr\}F_1''\biggl]\notag\\
&&+(\omega_1-\omega_2)^2\biggl[\bigr\{2\bigl\{J_{13}J_{14}F_2+(J_{13}J_{24}+J_{14}J_{23})F_2'\bigr\}(\omega_1\omega_2)^{1/2}\notag\\
&&-2\bigr\{\frac{J_{21}J_{22}}{(\omega_1\omega_2)^{1/2}}-J_{23}J_{24}(\omega_1\omega_2)^{1/2}\bigr\}F_2''\biggr]-(\omega_1-\omega_2)\biggl[\bigl\{J_{13}J_{22}(\frac{\omega_1}{\omega_2})^{1/2}\notag\\
&&-J_{14}J_{21}(\frac{\omega_2}{\omega_1})^{1/2}\bigr\}F_3'+2\bigl\{J_{21}J_{22}(\frac{\omega_2}{\omega_1})^{1/2}+J_{22}J_{23}(\frac{\omega_1}{\omega_2})^{1/2}\bigr\}F_3''\biggl]\notag\\
&&-\biggl[\bigr\{2\bigl\{J_{13}J_{14}F_4+(J_{13}J_{24}+J_{14}J_{23})F_4'\bigr\}(\omega_1\omega_2)^{1/2}-2\bigr\{\frac{J_{21}J_{22}}{(\omega_1\omega_2)^{1/2}}-J_{23}J_{24}(\omega_1\omega_2)^{1/2}\bigr\}F_4''\biggr]\Biggr\}\notag
\end{eqnarray}
\begin{eqnarray}
&r_7=&\frac{1}{3\omega_1^2(4\omega_1^2-\omega_2^2)} \Biggl\{
8\omega_1^3\biggl[J_{13}(J_{13}F_1+J_{23}F_1')\omega_1-\biggl(\frac{J_{21}^2}{\omega_1}-J_{23}^2\omega_1\biggr)F_1''\biggr]\\&&-2\omega_1\biggl[\omega_1J_{13}(J_{13}F_3+J_{23}F_3')-\biggl(\frac{J_{21}^2}{\omega_1}-J_{23}^2\omega_1\biggr)F_3''\biggr]\notag\\&&-4\omega_1^2J_{21}(J_{13}F_2+J_{23}F_2'')\omega_1+J_{21}(J_{13}F_4'+2J_{23}F_4'')\Biggr\}\notag
\end{eqnarray}
\begin{eqnarray}
&r_8&=\frac{-1}{3\omega_2^2(4\omega_2^2-\omega_1^2)} \Biggl\{
8\omega_2^3\biggl[J_{14}(J_{14}F_1+J_{24}F_1')\omega_2-\biggl(\frac{J_{22}^2}{\omega_2}-J_{24}^2\omega_2\biggr)F_1''\biggr]\\
&&+4\omega_2^2J_{22}(J_{14}F_2+2J_{24}F_2'')\omega_2-2\omega_2\biggl[\omega_2J_{14}(J_{14}F_3+J_{24}F_3')\notag\\ &&-\biggl(\frac{J_{22}^2}{\omega_2}-J_{24}^2\omega_2\biggr)F_3''\biggr]-J_{22}(J_{14}F_4'+2J_{24}F_4'')\Biggr\}\notag
\end{eqnarray}
\begin{eqnarray}
&r_9=&\frac{1}{\omega_1\omega_2(2\omega_1+\omega_2)(\omega_1+2\omega_2)}
\Biggl\{(\omega_1+\omega_2)^3\biggl[\bigr\{2J_{13}J_{14}F_1\\
&&+(J_{13}J_{24}+J_{14}J_{23})F_1'\bigr\}(\omega_1\omega_2)^{1/2}+2\bigr\{\frac{J_{21}J_{22}}{(\omega_1\omega_2)^{1/2}}+J_{23}J_{24}(\omega_1\omega_2)^{1/2}\bigr\}F_1''\biggr]\notag\\
&&-(\omega_1+\omega_2)^2\biggl[\bigl\{J_{13}J_{22}(\frac{\omega_1}{\omega_2})^{1/2}-J_{14}J_{21}(\frac{\omega_2}{\omega_1})^{1/2}\bigr\}F_2'\notag\\
&&-2\bigl\{J_{21}J_{24}(\frac{\omega_2}{\omega_1})^{1/2}-J_{22}J_{23}(\frac{\omega_1}{\omega_2})^{1/2}\bigr\}F_2''\biggl]
\notag\\
&&-(\omega_1+\omega_2)\biggl[\bigr\{2\bigl\{J_{13}J_{14}F_3+(J_{13}J_{24}+J_{14}J_{23})F_3'\bigr\}(\omega_1\omega_2)^{1/2}\notag\\
&&+2\bigr\{\frac{J_{21}J_{22}}{(\omega_1\omega_2)^{1/2}}+J_{23}J_{24}(\omega_1\omega_2)^{1/2}\bigr\}F_3''\biggr]-\biggl[\bigl\{J_{13}J_{22}(\frac{\omega_1}{\omega_2})^{1/2}\notag\\
&&-J_{14}J_{21}(\frac{\omega_2}{\omega_1})^{1/2}\bigr\}F_4'
-2\bigl\{J_{21}J_{24}(\frac{\omega_2}{\omega_1})^{1/2}-J_{22}J_{23}(\frac{\omega_1}{\omega_2})^{1/2}\bigr\}F_4''\biggl]
\Biggr\}\notag
\end{eqnarray}
\begin{eqnarray}
 & r_{10}=&\frac{1}{\omega_1\omega_2(2\omega_1-\omega_2)(2\omega_2-\omega_1)}
\Biggl\{(\omega_1-\omega_2)^3\biggl[\bigr\{2J_{13}J_{14}F_1\notag\\
&&+(J_{13}J_{24}+J_{14}J_{23})F_1'\bigr\}(\omega_1\omega_2)^{1/2}-2\bigr\{\frac{J_{21}J_{22}}{(\omega_1\omega_2)^{1/2}}-J_{23}J_{24}(\omega_1\omega_2)^{1/2}\bigr\}F_1''\biggr]\\
&&-(\omega_1-\omega_2)^2\biggl[\bigl\{J_{13}J_{22}(\frac{\omega_1}{\omega_2})^{1/2}-J_{14}J_{21}(\frac{\omega_2}{\omega_1})^{1/2}\bigr\}F_2'+2\bigl\{J_{21}J_{24}(\frac{\omega_2}{\omega_1})^{1/2}+J_{22}J_{23}(\frac{\omega_1}{\omega_2})^{1/2}\bigr\}F_2''\biggl]
\notag\\
&&-(\omega_1-\omega_2)\biggl[\bigr\{2\bigl\{J_{13}J_{14}F_3+(J_{13}J_{24}+J_{14}J_{23})F_3'\bigr\}(\omega_1\omega_2)^{1/2}-2\bigr\{\frac{J_{21}J_{22}}{(\omega_1\omega_2)^{1/2}}-J_{23}J_{24}(\omega_1\omega_2)^{1/2}\bigr\}F_3''\biggr]\notag\\
&&+\biggl[\bigl\{J_{13}J_{22}(\frac{\omega_1}{\omega_2})^{1/2}-J_{14}J_{21}(\frac{\omega_2}{\omega_1})^{1/2}\bigr\}F_4'+2\bigl\{J_{21}J_{24}(\frac{\omega_2}{\omega_1})^{1/2}-J_{22}J_{23}(\frac{\omega_1}{\omega_2})^{1/2}\bigr\}F_4''\biggl]
\Biggr\}\notag
\end{eqnarray}

We can write expressions of $s_i$ with the help of  $r_i$  replacing   $F_i$ by $G_i$, $F_i'$ by $G_i'$  and $F_i''$ by $G_i''$,$(i=1,2,3,4)$ where

\begin{eqnarray}
 F_1&=&\frac{-\zwe}{6}\\
 F_2&=&\frac{3}{32}\biggl[\frac{16}{3}\epsilon+6A_2-\frac{979}{18}\zae+\frac{(143+9\gamma)}{6\sqrt{3}}nW_1+\frac{(555+376\gamma)}{27\sqrt{3}}\zwe\\&&+\gamma\left\{14+\frac{4\epsilon}{3}+25A_2-\frac{1507 }{18}\zae-\frac{(215+29\gamma)}{6\sqrt{3}}nW_1
-\frac{2(1174+169\gamma)}{27\sqrt{3}}\zwe\right\}\biggr]\notag\\
F_3&=&\frac{3\sqrt{3}}{16}\biggl[14-\frac{16}{3}\epsilon+\frac{23A_2}{2}-\frac{104}{9}\zae+\frac{115(1+\gamma)}{18\sqrt{3}}nW_1-\frac{2(439-68\gamma)}{27\sqrt{3}}\zwe\\&&+\gamma\left\{\frac{32\epsilon}{3}+40A_2-\frac{310 }{9}\zae+\frac{(511+53\gamma)}{6\sqrt{3}}nW_1
-\frac{(2519-249\gamma)}{27\sqrt{3}}\zwe\right\}\biggr]\notag\\
F_4&=&\frac{-3}{256}\biggl[364+420A_2-\frac{17801A_2}{9}\zae+\frac{(2821+189\gamma)}{3\sqrt{3}}nW_1-\frac{(23077+9592\gamma)}{27\sqrt{3}}\zwe\\&&+28\gamma\left\{23+\frac{100\epsilon}{21}+\frac{849A_2}{14}+\frac{59 }{7}\zae-\frac{(125+38\gamma)}{6\sqrt{3}}nW_1
-\frac{(87613-213\gamma)}{27\sqrt{3}}\zwe\right\}\biggr]\notag
\end{eqnarray}
\begin{eqnarray}
 F_1'&=&\frac{\zwe}{3\sqrt{3}}\\
 F_2'&=&\frac{3\sqrt{3}}{16}\biggl[14-\frac{16}{3}\epsilon+A_2-\frac{1367}{18}\zae+\frac{115(1+\gamma)}{18\sqrt{3}}nW_1-\frac{(863-136\gamma)}{27\sqrt{3}}\zwe\\&&+\gamma\left\{\frac{32\epsilon}{3}+40A_2-\frac{382}{9}\zae+\frac{(511+53\gamma)}{6\sqrt{3}}nW_1
-\frac{(2519-24\gamma)}{27\sqrt{3}}\zwe\right\}\biggr]\notag\\
F_3'&=&\frac{-9}{8}\biggl[\frac{8}{3}\epsilon+\frac{203A_2}{6}-\frac{721}{54}\zae-\frac{(105+15\gamma)}{18\sqrt{3}}nW_1-\frac{(319-114\gamma)}{81\sqrt{3}}\zwe\\&&+\gamma\left\{2-\frac{4\epsilon}{9}-\frac{173A_2}{6}-\frac{781}{9}\zae+\frac{(197+23\gamma)}{18\sqrt{3}}nW_1
-\frac{(265-32\gamma)}{81\sqrt{3}}\zwe\right\}\biggr]\notag\\
F_4'&=&\frac{-3\sqrt{3}}{16}\biggl[392-\frac{532\epsilon}{3}+\frac{1918A_2}{3}-\frac{28582A_2}{9}\zae+\frac{(203+1211\gamma)}{9\sqrt{3}}nW_1+\frac{(949+4378\gamma)}{27\sqrt{3}}\zwe\\&&+28\gamma\left\{\frac{108\epsilon}{7}+\frac{4037A_2}{84}-\frac{611}{21}\zae+\frac{(8397+919\gamma)}{84\sqrt{3}}nW_1
-\frac{(92266-1869\gamma)}{27\sqrt{3}}\zwe\right\}\biggr]\notag
\end{eqnarray}
\begin{eqnarray}
 F_1''&=&\frac{\zwe}{6}\\
 F_2''&=&\frac{-9}{32}\biggl[\frac{8}{3}\epsilon+\frac{203A_2}{6}-\frac{625}{54}\zae-\frac{(105+15\gamma)}{18\sqrt{3}}nW_1-\frac{(307-114\gamma)}{81\sqrt{3}}\zwe\\&&+\gamma\left\{2-\frac{4\epsilon}{9}+\frac{55A_2}{2}-\frac{797}{54}\zae+\frac{(197+23\gamma)}{18\sqrt{3}}nW_1
-\frac{(211-32\gamma)}{81\sqrt{3}}\zwe\right\}\biggr]\notag\\
F_3''&=&\frac{-9\sqrt{3}}{16}\biggl[2-\frac{8}{3}\epsilon+\frac{55A_2}{6}-\frac{134}{3}\zae-\frac{(37+\gamma)}{18\sqrt{3}}nW_1-\frac{(93+226\gamma)}{81\sqrt{3}}\zwe\\&&+\gamma\left\{4\epsilon+\frac{169 }{27}\zae+\frac{(241+45\gamma)}{18\sqrt{3}}nW_1
-\frac{(1558-126\gamma)}{81\sqrt{3}}\zwe\right\}\biggr]\notag\\
F_4''&=&\frac{9}{256}\biggl[\frac{212}{3}\epsilon+\frac{2950A_2}{3}-\frac{1370A_2}{27}\zae-\frac{(771+237\gamma)}{9\sqrt{3}}nW_1-\frac{2(1907-984\gamma)}{81\sqrt{3}}\zwe\\&&+28\gamma\left\{\frac{11}{7}+\frac{4\epsilon}{9}-\frac{152A_2}{7}-\frac{36965}{504}\zae+\frac{(2569+277\gamma)}{252\sqrt{3}}nW_1
+\frac{(22603+4396\gamma)}{1134\sqrt{3}}\zwe\right\}\biggr]\notag
\end{eqnarray}
\begin{eqnarray}
 G_1&=&\frac{-\zwe}{6}\\
 G_2&=&\frac{3}{32}\biggl[14-\frac{16}{3}\epsilon+A_2-\frac{1367}{18}\zae+\frac{115(1+\gamma)}{18\sqrt{3}}nW_1-\frac{(863-136\gamma)}{27\sqrt{3}}\zwe\\&&+\gamma\left\{\frac{32\epsilon}{3}+40A_2-\frac{382 }{9}\zae+\frac{(511+53\gamma)}{6\sqrt{3}}nW_1-\frac{(2519-24\gamma)}{27\sqrt{3}}\zwe\right\}\biggr]\notag\\ 
G_3&=&\frac{3\sqrt{3}}{16}\biggl[\frac{16}{3}\epsilon+6A_2-\frac{907A_2}{18}\zae+\frac{(143+9\gamma)}{6\sqrt{3}}nW_1+\frac{(477+403\gamma)}{27\sqrt{3}}\zwe\\&&+\gamma\left\{14+\frac{4\epsilon}{3}+\frac{71A_2}{2}-\frac{1489}{18}\zae-\frac{(215+29\gamma)}{6\sqrt{3}}nW_1
-\frac{2(1174+169\gamma)}{27\sqrt{3}}\zwe\right\}\biggr]\notag\\
G_4&=&\frac{3\sqrt{3}}{256}\biggl[84+52\epsilon+212A_2-267\zae+\frac{2(299+61\gamma)}{3\sqrt{3}}nW_1-\frac{(14854+225\gamma)}{27\sqrt{3}}\zwe\\&&+\gamma\left\{32\epsilon+156A_2+649\zae-\frac{(562+8\gamma)}{3\sqrt{3}}nW_1+\frac{(13285+5169\gamma)}{27\sqrt{3}}\zwe\right\}\biggr]\notag
 \end{eqnarray}
\begin{eqnarray}
 G_1'&=&\frac{-\zwe}{\sqrt{3}}\\
 G_2'&=&\frac{9}{16}\biggl[\frac{8}{3}\epsilon+\frac{203A_2}{6}-\frac{625}{54}\zae-\frac{(105+15\gamma)}{18\sqrt{3}}nW_1-\frac{(307-114\gamma)}{81\sqrt{3}}\zwe\\&&-\gamma\left\{2-\frac{4\epsilon}{9}-\frac{55A_2}{2}-\frac{797}{54}\zae+\frac{(197+23\gamma)}{18\sqrt{3}}nW_1
-\frac{(211-32\gamma)}{81\sqrt{3}}\zwe\right\}\biggr]\notag\\
G_3'&=&\frac{3\sqrt{3}}{8}\biggl[14-\frac{16}{3}\epsilon+\frac{65A_2}{6}-\frac{1439}{18}\zae+\frac{115(1+\gamma)}{18\sqrt{3}}nW_1-\frac{(941-118\gamma)}{27\sqrt{3}}\zwe\\&&+\gamma\left\{\frac{32\epsilon}{3}-40A_2-\frac{310}{9}\zae+\frac{(511+53\gamma)}{6\sqrt{3}}nW_1
-\frac{(251-24\gamma)}{27\sqrt{3}}\zwe\right\}\biggr]\notag\\
G_4'&=&\frac{-9}{128}\biggl[12\epsilon-287A_2+\frac{847A_2}{9}\zae-\frac{2(28+\gamma)}{\sqrt{3}}nW_1-\frac{4(2210-69\gamma)}{27\sqrt{3}}\zwe\\&&-\gamma\left\{96+\frac{152\epsilon}{3}+135A_2-\frac{2320}{9}\zae+\frac{(497-123\gamma)}{3\sqrt{3}}nW_1
-\frac{4(17697+32\gamma)}{27\sqrt{3}}\zwe\right\}\biggr]\notag
\end{eqnarray}
\begin{eqnarray}
 G_1''&=&\frac{-\zwe}{6}\\
 G_2''&=&\frac{9\sqrt{3}}{32}\biggl[2-\frac{8}{3}\epsilon+\frac{23A_2}{3}-44\zae-\frac{(37+\gamma)}{18\sqrt{3}}nW_1-\frac{(123+349\gamma)}{3\sqrt{3}}\zwe\\&&+\gamma\left\{4\epsilon+\frac{88A_2}{27}+\frac{(421+45\gamma)}{18\sqrt{3}}nW_1
-\frac{(1558-126\gamma)}{81\sqrt{3}}\zwe\right\}\biggr]\notag\\
G_3''&=&\frac{-9}{16}\biggl[\frac{8}{9}\epsilon+\frac{203A_2}{6}-\frac{589}{54}\zae-\frac{5(51+2\gamma)}{18\sqrt{3}}nW_1-\frac{(349-282\gamma)}{81\sqrt{3}}\zwe\\&&+\gamma\left\{2-\frac{4\epsilon}{9}-26A_2-\frac{412}{27}\zae+\frac{(197+23\gamma)}{18\sqrt{3}}nW_1
-\frac{(211-32\gamma)}{81\sqrt{3}}\zwe\right\}\biggr]\notag\\
G_4''&=&\frac{-9\sqrt{3}}{256}\biggl[12+\frac{20}{3}\epsilon+76A_2-\frac{350A_2}{3}\zae+\frac{(32\gamma)}{3\sqrt{3}}nW_1-\frac{2(1529+450\gamma)}{27\sqrt{3}}\zwe\\\label{eq:g4dd}&&+\gamma\left\{8\epsilon-\frac{749A_2}{3}+\frac{808}{9}\zae-\frac{(109-40\gamma)}{3\sqrt{3}}nW_1
+\frac{(35-1269\gamma)}{27\sqrt{3}}\zwe\right\}\biggr]\notag
\end{eqnarray}

\section{Conclusion}
Using transformation $x=B_1^{\az}+B_2^{\az}$ and $y=B_1^{\za}+B_2^{\za}$ in Eq.(~\ref{eq:h3}) the third order part $H_3$ of the Hamiltonian in $I_1^{1/2}I_2^{1/2}$ is  of the form
\begin{equation}\label{eq:H3} H_3=A_{3,0}I_1^{3/2}+A_{2,1}I_1I_2^{1/2}+A_{1,2}I_1^{1/2})I_2+A_{0,3}I_2^{3/2}
\end{equation}
We can verify that in Eq.(~\ref{eq:H3}) $A_{3,0}$ vanishes independently as in Deprit and Deprit Barhtolom\'{e}(1967). Similarly the other coefficients $A_{2,1},A_{1,2},A_{0,3}$ are also found to be zero independently. Hence the third order part $H_3$ of the Hamiltonian in $I_1^{1/2}I_2^{1/2}$ is zero.

\thanks{\bf Acknowledgment:}{ We are thankful to D.S.T. Government of India, New Delhi for sanctioning a project DST/MS/140/2K dated 02/01/2004 on this topic. We are also thankful to IUCAA Pune for providing  financial assistance for visiting library and  computer facility.}

\end{document}